\documentclass[a4paper,12pt]{article}

\newcommand{\Z}{\mathbb{Z}}

\usepackage{amsmath}
\usepackage{amsthm}
\usepackage{amsfonts}
\usepackage{amssymb}
\usepackage{enumitem}

\title{\vspace{-40mm} \normalsize\textsc{Research Proposal on Percolation and Random Walks}\vspace{-5mm}}
\author{\normalsize{\text{A. Tzioufas}}}
\date{}
\begin{document}
\maketitle
\vspace{2mm}

\emph{Percolation theory} seeks to understand phenomena raging from epidemics and forest fires to the distribution of matter in the universe. Percolation refers to the canonical model for describing the flow of oil in disordered porous media. Percolation on integer lattices is the simplest model to undergo a \textit{phase transition}, that is, increasing local connectivity rules results in passing from connected clusters of exponentially small size, to the emergence of a (unique) \textit{infinite connected cluster} (Grimmett \cite{GRI, G} and Kesten \cite{K}). It is for this reason that the percolation cluster constitutes the canonical paradigm of an infinite random graph. The \emph{Contact Process} (Liggett \cite{L, L99}) and its discrete time analogue, viz.\ \emph{oriented percolation}, are the canonical percolation models incorporating a notion of: either time, when modeling an epidemic or, when modeling flow in a random medium, gravity.

\textit{Random walks} constitute the archetypical mathematical formalization of disordered motion resulting from successive random increments, and are a central topic in probability theory since the beginnings of the subject, with the investigations of De Moivre, Laplace, Bernoulli, Pascal, among others. The main references on the topic, whose everlasting influence cannot be overestimated, are those of Feller [Chpts.\ III and XIV, \cite{F68}] and Spitzer \cite{S76}, whereas, from the plethora of more recent ones, we refer to Lawler \cite{L10}, Lawler and Limic \cite{LL10} due to being closer to our scope, and also to the amusing introduction of Doyle and Snell \cite{Doy}. We will exploit connections of random walks with the so-called \textit{Brownian motion} and with \textit{electrical networks}. 

Brownian motion (or pedesis) refers to the phenomenon of chaotic displacements of small particles suspended in a liquid or gas resulting from collisions with the molecules of the medium, the existence of which was empirically confirmed by botanist Robert Brown in 1827, although known earlier references to it date back to scientific poem of poet and philosopher Lucretius "On the Nature of Things" (c. 60 BC). Its precise mathematical explanation, offered by Albert Einstein in one of his Annus mirabilis papers, 1905, served as definitive confirmation that atoms and molecules, as theorized by the so-called Atomists since ancient times, actually exist.  Although there exist several mathematical models for the phenomenon, \textit{Brownian motion} (or the \textit{Wiener process}) is one of the most important ones in the theory of random processes. It arises as the scaling limit when speeding up the walk and taking the lattice spacing go to zero appropriately. 

Furthermore, we shall exploit the intimate, deep connection amongst random walks and the classical physics subject of \textit{electrical networks}. This connection has had profound consequences and, indeed, the interplay between the two has been proven mutually beneficial. 
It allows for probabilists to draw on a large body of ways of thought and well established methods from the physics literature, most prominently involving considerations of energy as the Thomson and Dirichlet principles. Whereas it also provides with insights and interpretations to quantities and physical laws pertaining to electrical networks, as for example that of effective resistance and its monotonicity law (cf.\ Doyle and Snell \cite{Doy}).

Part 1 of this proposal regards work on the critical contact process. 
Part 2 regards work on the random walk and Brownian motion. 
Part 3 regards work on a simple stochastic neural networks model inspired from percolation.
Part 4 lies in the interface of percolation and random walks, in connection with electrical networks. 
In what follows we elaborate on each of the parts of the project.  



\begin{enumerate}[leftmargin=0cm,itemindent=.5cm,labelwidth=\itemindent,labelsep=0cm,align=left]

	\item{
	
	Since the foundational work of Bezuidenhout and Grimmett \cite{BG0} it is known that the contact process on integer lattices possesses no more than one non-trivial invariant measure. This is a consequence of that the process from the fully occupied configuration at criticality converges to the empty configuration in a weak sense. In this work we show that, nevertheless, each finite subset is at a fully occupied state for arbitrarily large times with probability one. This is a manifestation of that weak convergence results fail to capture phenomena that almost sure ones do. The method of proof we use for extending this consequence of the non-negative asymptotic speed result in dimensions higher than one relies on techniques developed in Bezuidenhout and Grimmett \cite{BG1} for establishing there that the critical exponent associated to the Lebesgue integral of the occupied region in the graphical representation of the process assumes value which is greater than or equal than that predicted by mean field theory (2, in the logarithmic sense). 


We investigate upon extensions and consequences of the aforementioned result. As a first consequence we show that a certain interacting particle system, the idle contact process at criticality, exhibits strong survival with positive probability (w.p.p.). In this process, at the outset, an idle particle is placed on each site of $\Z^{d}$ other than the origin, which is inhabited by an excited particle instead. Idle particles are set to diffuse their descendancy according to a contact process at criticality, if ever excited. Once excited an idle particle immediately initiates reproduction, as well as attempts to excite its neighbors at a certain rate before dying. Progeny of different particles evolves otherwise independent of one another. The proof of that this process exhibits strong survival w.p.p.\ is a consequence of an extension of the result mentioned in a) in the case that the initial configuration is that of an infinite (supercritical) site percolation cluster and by means of offering a new basic coupling for independent contact processes. To show the former mentioned extension we rely on that an infinite percolation cluster does not fail to include an infinite subset of any fixed on the outset infinite collection of points, a consequence of its uniqueness by basic ergodicity properties. Two other corollaries regarding weakened initial conditions under which the result in a) remains valid for configurations that contain an infinite number of infinite, finite width strips of vacant sites, and for configurations distributed according to the invariant measure of a highly supercritical contact process are shown as byproducts.}



 

\item{This part is concerned with revisiting a problem of quite long mathematical history. The so-called \textit{absorption problem} regards asymptotics of the random times the module of the walk attains new maxima values, i.e.\ the times the walk exits symmetric intervals about its origin\footnote{the problem's nomenclature derives from the equivalent perspective of the walk restrained on intervals with absorbing endpoints} (Theorem 2.13 Revesz \cite{R}). As an alternative to the various available interesting analytical methods, we show a new approach that is elementary in its entirety. Our approach relies on Laplace transformations apparatus and leans on the so-called continuity theorems, commonly attributed to P.\ Levy.  We find this approach of interest in its own right due to its overall rudimentary and simple nature. An additional interesting feature of this approach is that it yields as byproduct a connection with first-passage times, and the stable law of order 1/2 limit for first-passages times is retrieved. 

We exploit this result to show an elementary proof to that symmetric planar random walks exit times from spheres and partial maxima values, under appropriate rescaling, possess associated asymptotic distributions, which are identified in explicit form and associated to functionals of standard planar Brownian motion. A link to the solution to associated boundary value problems of the simplest type is also presented.}


\item{
We purpose understanding functions associated to observed central nervous system phenomena
with aid of a new, conceptually simplifying class of models. We will argue that these models, despite differentiating to a great respect from existing ones, possess characteristics qualifying them as realistically befitting alternatives.

Cognitive science refers to the study of processes of the mind, more in particular, the development of behavior and intelligence, focusing especially on how information is represented, processed, and transformed in nervous systems (humans and animals) or machines (computers). It thus encompasses various scientific fields as, for instance, neuroscience, linguistics, artificial intelligence, or anthropology. Cells, the so-called "building blocks of life" due to being the smallest unit of life which replicates independently,
are the basic functional, structural, biological units of all known living organisms. The \textit{central nervous system} (CNS), viz.\ the human brain and spinal cord, is a prominent example of a biological neural network. {\it Neural networks} is a generic term referring to the theoretical framework of models that aim at understanding, explaining, and predicting the behavior of such complex systems. 
A \textit{neuron} is an electrically excitable cell that processes and transmits information through electrical and chemical signals with other cells via \textit{synapses}; in our model below, although no differentiation between electrical and chemical ones will be made, bidirectional synapses will be incorporated\footnote{it is known however that only the former (commonly) possess this property.}.  
A key feature such models seek to understand is the synchronized activity of large numbers of neurons resulting to oscillations of neural ensembles\footnote{fact that can be empirically observed through experimental measurements known as the \textit{electroencephalogram}} (aka \textit{spike trains}). An action potential is an \textit{instantaneous} event in which the electrical membrane potential of a cell rapidly rises and falls. The temporal sequence of action potentials generated by a neuron is called its \textit{spike train}.


Here, we will introduce a neural network process with long-range dynamic synaptic connections. This model incorporates ideas from both standard neural network models of the simplest type (for instance, the Hopfield model cf.\ Rojas \cite{R96}, see also Cessac \cite{C11}), as well as dynamical long-range percolation processes (cf.\ Grimmett \cite{G}; H\"aggstr\"om, Peres and Steif \cite{HPS97}). One additional innovative aspect of our model will be that occurrence of newly observed spikes will result from spike trains arriving from infinity, i.e.\ entire 
trails of electrical activity among consecutively connected neurons arriving at individual neurons from infinity.
Long-range percolation on the integers is a probability model in which any two integer sites of some countable space are connected with some probability inversely proportional to their distance to the power of parameter value, whereas nearest connections among neighboring sites occur with some distinct probability parameter. These connections will be interpreted in our model as representing synaptic connections. We aim at identifying criteria for the parameters associated to the renewable synaptic connections that guarantee perpetually non-trivial activity of the process, in fact, we will identify so-called {\it phase transition} phenomena depending on their values. 

The adult human brain is estimated to contain  $\sim 10^{11}$ neurons which is sufficiently large to justify consideration of infinite dimensional neural networks for modeling purposes. In addition however, it contains a huge number (from $10^{14}$ to $5 \times 10^{14}$) of synapses\footnote{every cubic millimeter of cerebral cortex contains roughly $10^{9}$ of them} (cf.\ with e.g.\ Drachman \cite{D05}), whereas recall that, the speed of electromagnetic waves is  
known to be of the magnitude of $\sim 3 \times 10^{8}$ m/sec, which is large when compared to the size of the brain of a mammal. 
These facts respectively explain our motivation for introducing long (unbounded) range synaptic connections, in addition to nearest neighbors ones and, furthermore, the conceptualizing of spike trains occurring according to synaptic activity arriving from infinity, equally well.  In addition, synaptic connections which evolve dynamically in random are introduced, i.e.\ are renewable according to an underlying random (percolation) process. Nonetheless, the evolution of our simple process is Markovian. For related non-Markovian probability models, {\it viz.} the so-called chains with memory of variable length, cf.\ with Fern\'andez, Ferrari, Galves \cite{FFG01}, as well as with Galves and L\"ocherbach \cite{GL13} for models specialized to neuronal networks.  


We will introduce an interacting particle system in dynamic, long-range percolation random environment process, which encompasses basic features and fundamental principles exhibited by biological neural networks. We will identify regions of the parameters measuring the strength of neuron connections where the system exhibits non-trivial activity. Furthermore, we aim to provide with stationarity results for the process, and with bounds on the interarrival so-called spiking times as functions of the values of these parameters. We expect our analysis to leverage from an ensemble of celebrated results on (static, non-oriented) long-range percolation due to Schulman \cite{S83}, Kalikow (cf.\ with Kalikow and Weiss \cite{KW88}),  Newman and Schulman \cite{NS86}, and Aizenman and Newman \cite{AN86}. 

}

\item{We consider the Simple Random Walk on an infinite percolation cluster model, \textit{viz.} the so-called \textit{ant in a labyrinth}.  
Results concerning the scaling limit of this process have been of central interest in the probability literature for years (cf.\ with the recent survey of Biskup \cite{B}). The problem under consideration, originally proposed in the monograph of Kesten \cite{K}, concerns the limiting behavior and asymptotics of the effective resistance of the electrical network associated with the cluster of supercritical bond percolation. By elementary considerations involving Kirchhoff's and Ohm's laws the effective resistance and voltages on the vertices of the network correspond to determining certain passage times for the random walk on this cluster and the construction of the harmonic function respectively. We investigate open problems 12 and 13 from Kesten \cite{K}. The conjectures effectively concern extending  to the case of a the walk on the random environment of an infinite percolation cluster the well-known asymptotic scaling result for escape probabilities of the Random Walk (see for instance Proposition 2.16 in \cite{LP}, or Lemma 22.1 in \cite{R}). To resolve these problems, we point out to known asymptotics regarding the geometry of the cluster, in particular, regarding the boundary to volume ratio, see Theorem 8.99 in Grimmett \cite{GRI}. The reason we believe this ratio hints the right scaling for the quenched escape probabilities correctly is that it seems to represent (intuitively) the quote of the escape paths remaining open (volume) to those blocked (boundary). }

\end{enumerate}



\begin{thebibliography}{99}
\footnotesize



\bibitem{AN86}  \textsc{Aizenman, M., and Newman, C. M.} Discontinuity of the percolation density in one dimensional $1/| x - y|^{2}$ percolation models. \textit{Communications in Mathematical Physics}  \textbf{107.4} (1986): 611-647.


\bibitem{BG0} { \sc Bezuidenhout, C.} and { \sc Grimmett, G.} (1990). The critical contact process dies out. \textit{Ann. Probab.} \textbf{18} 1462--1482.

\bibitem{BG1} { \sc Bezuidenhout, C. and Grimmett, G.} (1991). Exponential decay for subcritical contact and percolation processes.  \textit{Ann. Probab.} \textbf{19} 984-1009 .

\bibitem{B} {\sc Biskup, M.} (2011). Recent progress on the Random Conductance Model. \textit{Prob. Surveys.}\textbf{ 8} 

\bibitem{BP} {\sc Billingsley, P.} (2008). Probability and measure. John Wiley \& Sons.


\bibitem{C11} \textsc{Cessac, B.} A discrete time neural network model with spiking neurons: II: Dynamics
with noise. \textit{Journal of Mathematical Biology}, \textbf{62} :863–900, 2011.

\bibitem{D05}   \textsc{Drachman, D.} (2005). Do we have brain to spare?. {\it Neurology} {\bf 64 (12)}: 2004–5. doi:10.1212/01.WNL.0000166914.38327.BB. PMID 15985565.


\bibitem{D} \textsc{Durrett, R.} (1984). Oriented percolation in two dimensions. \textit{Ann. Probab.} 999-1040.


\bibitem{D91}  \textsc{Durrett, R.} (1991). The contact process, 1974-1989. Cornell University, Mathematical Sciences Institute.

\bibitem{Doy} {\sc Doyle, P.} and {\sc Snell,  L.} (1984). Random walks and electric networks. Vol. 22. Mathematical association of America.

\bibitem{F68}  \textsc{Feller, W.} {\em An introduction to probability theory and its applications.} Vol. 1. John Wiley \& Sons, 1968.


\bibitem{FFG01} \textsc{Fern\'andez, R., Ferrari, P.A., Galves, A.} Coupling, renewal and perfect simulation of chains of infinite order. \& 





\bibitem{GL13} \textsc{Galves, A., and L\"ocherbach, E.} Infinite systems of interacting chains with memory of variable length—a stochastic model for biological neural nets. \textit{Journal of Statistical Physics} \textbf{151.5} (2013): 896-921.


\bibitem{GRI} { \sc Grimmett, G.R.} (1999). {\em Percolation.} Springer, Berlin.


\bibitem{G}  {\sc Grimmett, G.R.} (2010). Probability on Graphs. Cambridge University Press.



\bibitem{HPS97} {\sc H\"aggstr\"om, O., Peres, Y. and Steif, J.}. Dynamical percolation. \textit{Annales de l'IHP Probabilités et Statistiques.} {\bf 33(4)}. 1997.










\bibitem{KW88} {\sc Kalikow, S., and B. Weiss.} "When are random graphs connected."  \textit{Israel journal of mathematics} {\bf 62.3 } (1988): 257-268.


\bibitem{K}  {\sc Kesten, H.} (1982). Percolation Theory for Mathematicians.






\bibitem{L10}  \textsc{Lawler, G.} {\em  Random walk and the heat equation.} Vol. 55. American Mathematical Soc., 2010.


\bibitem{LL10}  \textsc{Lawler, G., and Limic, V.} {\em Random walk: a modern introduction.} Vol. 123. Cambridge University Press, 2010.


\bibitem{L}  {\sc Liggett, T.} (1985). {\em Interacting particle systems.} Springer, New York. 
 
\bibitem{L99}  {\sc  Liggett, T.} (1999). {\em Stochastic Interacting Systems: Contact, Voter and Exclusion Processes.} Springer, New York.







\bibitem{LP}  {\sc Levin, A.},  {\sc Peres, W.,} and  {\sc Wilmer, E.} (2009). Markov Chains and Mixing Times 

\bibitem{NS86} \textsc{Newman, C. M., and Schulman, L. S.} One dimensional $1/| j - i|^{s}$ percolation models: The existence of a transition for $s \leq 2$. \textit{Communications in Mathematical Physics},  \textbf{104(4)}, (1986): 547-571.




\bibitem{R96}{ \sc Rojas, R.} {\it Neural Networks - A Systematic Introduction}. Springer-Verlag, Berlin, New-York, 1996.



\bibitem{R} {\sc Revesz, P.} (1990). Random Walk in Random and non-Random environment. 

\bibitem{S83} {\sc Schulman, L. S.} (1983). Long range percolation in one dimension. \textit{Journal of Physics A: Mathematical and General}, {\bf 16(17)}, L639.
ISO 690	


\bibitem{S76}  \textsc{Spitzer, F.} {\em Principles of random walk.} Vol. 34. Springer 1976.


\end{thebibliography}
\end{document}